\documentclass[11pt,a4paper]{article}
\usepackage{amssymb}
\usepackage{amsfonts}

\newtheorem{theorem}{Theorem}

\newtheorem{proposition}[theorem]{Proposition}

\begin{document}

\title{A characterization of the Khavinson-Shapiro conjecture via Fischer
operators}
\date{}
\author{Hermann Render \\
\noindent School of Mathematics and Statistics, \noindent University College
Dublin, \\
Belfield, Dublin 4, Ireland. Email: hermann.render@ucd.ie}
\maketitle

\begin{abstract}
The Khavinson-Shapiro conjecture states that ellipsoids are the only bounded
domains in euclidean space satisfying the following property (KS): the
solution of the Dirichlet problem for polynomial data is polynomial. In this
paper we show that property (KS) for a domain $\Omega $ is equivalent to the
surjectivity of a Fischer operator associated to the domain $\Omega .$
\end{abstract}

\section{Introduction}

\footnotetext[1]{%
2010 \textit{Mathematics Subject Classification } 31B05; 35J05. Keywords and
phrases: Dirichlet problem, harmonic extension, Khavinson-Shapiro conjecture.%
}

In the 19th century ellipsoidal harmonics have been used to prove that for
any polynomial $p$ of degree $\leq m$ there exists a harmonic polynomial $h$
of degree $\leq m$ such that $h\left( \xi \right) =p\left( \xi \right) $ for
all $\xi \in \partial E$ where $\partial E $ is the boundary of an \emph{%
ellipsoid} $E$ in the euclidean space $\mathbb{R}^{d}.$ It follows that an
ellipsoid satisfies the following property defined for an arbitrary open
subset $\Omega $ in $\mathbb{R}^{d}$:

\begin{itemize}
\item[(KS)] For any polynomial $p$ with real coefficients there exists a
harmonic polynomial $h$ with real coefficients such that $h\left( \xi
\right) =p\left( \xi \right) $ for all $\xi \in \partial \Omega .$
\end{itemize}

The Khavinson-Shapiro conjecture \cite{KhSh92} states that ellipsoids are
the only bounded domains $\Omega $ in $\mathbb{R}^{d}$ with property (KS).
Obviously a domain $\Omega $ has property (KS) if and only if the Dirichlet
problem for polynomial data (restricted to the boundary) have polynomial
solutions; for the Dirichlet problem we refer to \cite{ArGa01} and \cite%
{Gard93}. The Khavinson-Shapiro conjecture has been confirmed for large
classes of domains but it is still unproven in its full generality, and we
refer the interested reader to the expositions \cite{HaSh94}, \cite{Render}, 
\cite{Lund}, \cite{KhLu10}, \cite{LuRe10}, \cite{KhLu14} and for further
ramifications in \cite{KhSt} originating from the work \cite{PutSty}.

In this paper we want to characterize the property (KS) by using Fischer
operators. In our context we shall mean by a Fischer operator\footnote[2]{%
More generaly one can define a Fischer operator by $q\longmapsto P\left(
D\right) \left( \psi q\right) $ where $P\left( D\right) $ is a linear
partial differential operator with constant real coefficients. Fischer's
Theorem in \cite{Fischer} states that the Fischer operator is a bijection
whenever $\psi \left( x\right) $ is a homogeneous polynomial equal to the
polynomial $P\left( x\right) .$} an operator of the form 
\[
F_{\psi }\left( q\right) :=\Delta \left( \psi q\right) \hbox{ for }q\in 
\mathbb{R}\left[ x\right] 
\]%
where $\Delta $ is the Laplace operator $\frac{\partial ^{2}}{\partial
x_{1}^{2}}+....+\frac{\partial ^{2}}{\partial x_{d}^{2}}$ and $\psi $ is a
fixed element in $\mathbb{R}\left[ x\right] $, the set of all polynomials in 
$d$ variables with real coefficients. Fischer operators often allows
elementary and short proofs of mathematical statements which usually require
hard and deep analysis, see \cite{Shap89}, \cite{KhLu14}. For example, the
statement that ellipsoids have property (KS) can be proven in a few lines
using Fischer operators and elementary results in Linear Algebra, see \cite%
{Ba}, \cite{ChSi01}, \cite{AGV03}, and for further generalizations see \cite%
{Arm}, \cite{KhSh92}.

In order to formulate our main result we need some technical definitions.
The zero-set of a polynomial $f\in \mathbb{R}\left[ x\right] $ is denoted by 
$Z\left( f\right) =\left\{ x\in \mathbb{R}^{d}:f\left( x\right) =0\right\} .$
We say that a subset $Z$ of $\mathbb{R}^{d}$ is an \emph{admissible common
zero set} if there exist \emph{non-constant irreducible} polynomials $f,g\in 
\mathbb{R}\left[ x\right] $ such that (i) $f\neq \lambda g$ for all $\lambda
\in \mathbb{R}$ and (ii) $Z=Z\left( f\right) \cap Z\left( g\right) .$ For
dimension $d=2$ it is well known that an admissible common zero set is
finite, see \cite[p. 2]{Shaf}. For arbitrary dimension $d$ it is intuitively
clear that an admissible common zero set has "dimension" $\leq d-2$ at each
point.

We say that an open set $\Omega $ in $\mathbb{R}^{d}$ is \emph{admissible}
if for any $x\in \partial \Omega $, any open neighborhood $V$ of $x$ and for
any finite family of admissible common zero sets $Z_{1},...,Z_{r}$ the set 
\[
\left[ \partial \Omega \cap V\right] \diagdown \bigcup_{j=1}^{r}Z_{j} 
\]%
is non-empty. For dimension $d=2$ it is easy to see that an open set $\Omega 
$ is admissible if each point $x\in \partial \Omega $ is not isolated in $%
\partial \Omega $. For arbitrary dimension it seems to be difficult to
formulate a precise topological condition but it is intuitively clear that a
domain $\Omega $ is admissible if each point in the boundary $\partial
\Omega $ has a neighborhood of dimension $d-1.$

The following is now the main result of this paper:

\begin{theorem}
\label{ThmI1}Let $\Omega $ be an open admissible subset of $\mathbb{R}^{d}$.
Then property (KS) holds for $\Omega $ if and only if there exists a
non-constant polynomial $\psi \in \mathbb{R}\left[ x\right] $ such that (i) $%
\partial \Omega \subset \psi ^{-1}\left\{ 0\right\} $ and (ii) the Fischer
operator $F_{\psi }:$ $\mathbb{R}\left[ x\right] \rightarrow \mathbb{R}\left[
x\right] $ defined by $F_{\psi }\left( q\right) :=\Delta \left( \psi
q\right) $ for $q\in \mathbb{R}\left[ x\right] $ is surjective.
\end{theorem}

An immediate consequence of Theorem \ref{ThmI1} is that the
Khavinson-Shapiro conjecture is true for all admissible bounded domains if
the following purely algebraic conjecture of M. Chamberland and D. Siegel
formulated in \cite{ChSi01} is true:

\begin{itemize}
\item[(CS)] The surjectivity of the Fischer operator $F_{\psi }:$ $\mathbb{R}%
\left[ x\right] \rightarrow \mathbb{R}\left[ x\right] $ implies that the
degree of $\psi $ is $\leq 2.$
\end{itemize}

We refer to \cite{LuRe10} and \cite{Rend15} for more details on conjecture
(CS)\ and related results. It should be emphasized that the polynomial $\psi 
$ in conjecture (CS) has real coefficients. In \cite{Khav91} it is shown for
dimension $d=3$ that for any non-constant polynomial $\varphi \left(
z\right) =a_{0}+a_{1}z+\cdots +a_{n}z^{n}$ in the complex variable $z$ the
operator $F_{\psi }:\mathbb{C}\left[ x\right] \rightarrow \mathbb{C}\left[ x%
\right] $ defined by 
\[
F_{\psi }\left( q\right) =\Delta \left[ \left( x_{3}-\varphi \left(
x_{1}+ix_{2}\right) \right) ^{2}q\left( x_{1},x_{2},x_{3}\right) \right] 
\]%
for $q\in \mathbb{C}\left[ x\right] $ is surjective where $\mathbb{C}\left[ x%
\right] $ is the set of all polynomials with complex coefficients.

\section{Proof of Theorem 1}

From \cite{LuRe10} we cite the following result which is known as the \emph{%
Fischer decomposition} of a polynomial.

\begin{proposition}
\label{connect} Suppose $\psi$ is a polynomial. Then the operator $F_{\psi
}\left( q\right) :=\Delta\left( \psi q\right) $ is surjective if and only if
every polynomial $f$ can be decomposed as $f=\psi q_{f}+h_{f}$, where $q_{f}$
is a polynomial and $h_{f}$ is harmonic polynomial
\end{proposition}

\textbf{Proof of Theorem 1:} The sufficiency is easy: By assumption, there
exists $\psi \in \mathbb{R}\left[ x\right] $ with $\partial \Omega \subset
\psi ^{-1}\left\{ 0\right\} $ such that $F_{\psi }$ is surjective. According
to Proposition \ref{connect} there exists for each polynomial $f$ a
polynomial $q$ and a harmonic polynomial $u$ such that $f=\psi q+u$. For $%
\xi \in \partial \Omega \subset \psi ^{-1}\left\{ 0\right\} $ it follows
that $f\left( \xi \right) =u\left( \xi \right) $. Thus $u$ is a polynomial
solution to the Dirichlet problem of the domain $\Omega $ and property (KS)
is satisfied.

Now assume that (KS) holds. Then there exists a harmonic polynomial $u\in 
\mathbb{R}\left[ x\right] $ such that $\left\vert \xi \right\vert
^{2}=u\left( \xi \right) $ for all $\xi \in \partial \Omega .$ Define the
polynomial $Q\left( x\right) =\left\vert x\right\vert ^{2}-u\left( x\right)
. $ Then 
\begin{equation}
\partial \Omega \subset Q^{-1}\left( 0\right)  \label{eqinc}
\end{equation}%
and $Q$ is a non-constant polynomial of degree $\geq 2$ since $\Delta \left(
Q\right) \neq 0.$ We factorize $Q\left( x\right) $ in irreducible factors, so%
\[
Q\left( x\right) =\left\vert x\right\vert ^{2}-u\left( x\right)
=f_{1}^{m_{1}}...f_{r}^{m_{r}} 
\]%
where $f_{k}$ is not a scalar multiple of $f_{l}$ for $k\neq l,$ and $%
m_{k}\geq 1$ is the multiplicity of $f_{k}.$ It follows that 
\begin{equation}
\partial \Omega \subset \bigcup_{k=1}^{r}f_{k}^{-1}\left( 0\right) .
\label{Gcontain}
\end{equation}%
Then $Z\left( f_{k}\right) \cap Z\left( f_{l}\right) $ is an admissible
common zero set for $k\neq l.$ Let $Z$ be the (finite) union of the sets $%
Z\left( f_{k}\right) \cap Z\left( f_{l}\right) $ with $k\neq l.$ As $\Omega $
is admissible there exists $x\in \partial \Omega \setminus Z,$ and now (\ref%
{Gcontain}) implies that 
\[
I:=\left\{ k\in \left\{ 1,...r\right\} :\exists x\in \partial \Omega
\setminus Z\hbox{ with }f_{k}\left( x\right) =0\right\} 
\]%
is non-empty. Let us define 
\[
\psi :=\prod_{k\in I}f_{k}\left( x\right) . 
\]%
We want to show that $F_{\psi }$ is surjective. By Proposition \ref{connect}
it suffices to show that for any polynomial $f$ there is a harmonic
polynomial $u$ and a polynomial $q$ such that $f=\psi q+u.$ By property (KS)
there exists a harmonic polynomial $u$ such that $f\left( \xi \right)
=u\left( \xi \right) $ for all $\xi \in \partial \Omega .$ Define $g\left(
x\right) =f\left( x\right) -u\left( x\right) ,$ so $g$ vanishes on $\partial
\Omega .$ If we can show that each $f_{k}$ with $k\in I$ divides $g$ then,
by irreducibility of $f_{k}$ and the condition that $f_{k}\neq \lambda f_{j}$
for $k\neq j,$ we infer that $\psi =\prod\limits_{k\in I}f_{k}$ divides $g,$
say $g=\psi q$ for some polynomial $q.$ Then $f-u=g=q\psi $ and we are done.

Let $k\in I$ be fixed and $g$ as above. Then there exists $x\in \partial
\Omega \setminus Z$ with $f_{k}\left( x\right) =0.$ Then $f_{l}\left(
x\right) \neq 0$ for all $l\neq k$ since otherwise $x$ would be in $Z.$ By
continuity there is an open neighborhood $V$ of $x$ such that $f_{l}\left(
y\right) \neq 0$ for all $y\in V$ and $l\neq k.$ Then we conclude from (\ref%
{Gcontain}) that 
\begin{equation}
\partial \Omega \cap V\subset f_{k}^{-1}\left( 0\right) .  \label{eqinclus}
\end{equation}%
Let us write $g=g_{1}^{m_{1}}\cdots g_{s}^{m_{r}}$ where $g_{1},...,g_{s}$
are irreducible polynomials such that $g_{j}\neq \lambda g_{l}$ for $j\neq
l. $ If $f_{k}=\lambda g_{j}$ for some $j\in \left\{ 1,...,s\right\} $ we
see that $f_{k}$ divides $g.$ Assume that this is not the case and define $%
Z_{k}$ as the (finite) union of the admissible sets $Z\left( g_{j}\right)
\cap Z\left( f_{k}\right) $ for $j=1,...s.$ Since $\Omega $ is admissible
there exists $y\in \partial \Omega \cap V\diagdown \left( Z\cup Z_{k}\right) 
$. The inclusion (\ref{eqinclus}) shows that $f_{k}\left( y\right) =0,$ and
since $g$ vanishes on $\partial \Omega $ there exists $j\in \left\{
1,..,s\right\} $ such that $g_{j}\left( y\right) =0.$ Hence $y\in Z\left(
g_{j}\right) \cap Z\left( f_{k}\right) \subset Z_{k}.$ Now we obtain a
contradiction since $y\in \partial \Omega \cap V\diagdown \left( Z\cup
Z_{k}\right) .$

It remains to prove that $\partial \Omega $ is contained in $\psi
^{-1}\left( 0\right) .$ If $j\in \left\{ 1,...r\right\} \setminus I$ then it
follows from the definition of $I$ that $f_{j}\left( x\right) \neq 0$ for
all $x\in \partial \Omega \setminus Z.$ This fact and (\ref{Gcontain}) imply
that 
\begin{equation}
\partial \Omega \setminus Z\subset \bigcup_{k\in I}f_{k}^{-1}\left( 0\right)
=:F.  \label{eqII}
\end{equation}%
Let $x\in \partial \Omega .$ Since $\Omega $ is admissible there exists for
any ball $V$ with center $x$ and radius $1/m$ an element $x_{m}\in (\partial
\Omega \cap V)\diagdown Z.$ Then (\ref{eqII}) shows that $x_{m}\in F.$ Since 
$x_{m}$ converges to $x$ and $F$ is closed we infer that $x\in F.$ Thus $%
\partial \Omega \subset F.$

\end{document}